\input amstex
\mag=\magstep1
\documentstyle{amsppt}

\topmatter
\title A Base Point Free Theorem of Reid Type, II
\endtitle
\author Shigetaka Fukuda
\leftheadtext{S. Fukuda}
\rightheadtext{Base point free theorem}
\endauthor
\endtopmatter

\document

\flushpar
Faculty of Education, Gifu Shotoku Gakuen University, Yanaizu-cho, Hashima-gun, Gifu-prefecture, 501-6194, Japan

\flushpar
(e-mail: fukuda\@ha.shotoku.ac.jp)

\medpagebreak

\specialhead Introduction
\endspecialhead

\flushpar
This paper is a continuation of [Fu2].

We generally use the notation and terminology of [Utah].

Let $X$ be a normal, complete algebraic variety over {\bf C} and $(X,\Delta)$ a log variety which is log canonical.
We assume that $K_X+\Delta$ is a {\bf Q}-Cartier {\bf Q}-divisor.
Let $r$ be the smallest positive integer such that $r(K_X+\Delta)$ is Cartier ($r$ is called the {\it singularity index} of $(X,\Delta)$).

\definition{Definition (due to Reid [Rd])}
Let $\Theta = \sum_{i=1}^s \Theta_i$ be a reduced divisor with only simple normal crossings on an $n$-dimensional non-singular complete variety over {\bf C}.
We denote {\bf Strata} $(\Theta):= \{ \Gamma \vert 1 \leq k \leq n, 1 \leq i_1 < i_2 < \cdot \cdot \cdot <i_k \leq s, \quad \Gamma$ is an irreducible component of $\Theta_{i_1} \cap \Theta_{i_2} \cap \cdot \cdot \cdot \cap \Theta_{i_k} \not= \emptyset \} .$

Let $f:Y \to X$ be a log resolution of $(X,\Delta)$ such that $K_Y = f^*(K_X+\Delta)+\sum a_j E_j$ (where $a_j \ge -1$).
Let $L$ be a Cartier divisor on $X$.
$L$ is called {\it nef and log big} on $(X,\Delta)$ if $L$ is nef and big and $(L \vert _{f(\Gamma)})^{dim f(\Gamma)}>0$ for any member $\Gamma$ of {\bf Strata}$(\sum_{a_j=-1}E_j).$
\enddefinition

We note that the definition of the notion of "nef and log big" does not depend on the choice of the log resolution $f$ ([Fu1, Claim]).
In the case in which $(X,\Delta)$ is Kawamata log terminal (klt), if $L$ is nef and big, then $L$ is nef and log big on $(X,\Delta)$.

In [Rd], M. Reid gave the following statement:

\flushpar
{\it Let $L$ be a nef Cartier divisor such that $aL-(K_X+\Delta)$ is nef and log big on $(X,\Delta)$ for some $a \in$ {\bf N}.
Then Bs$\vert mL \vert = \emptyset$ for every $m \gg 0$.}

\flushpar
In the case in which $(X,\Delta)$ is klt, this is the standard Kawamata-Shokurov result (cf. [Ka84, Theorem 2.6]).
While in the case in which $(X,\Delta)$ is weakly Kawamata log terminal, the proof in [Fu2] of the statement needs the log minimal model program, which is still a conjecture in dimension $\ge 4$ (the assumption that $X$ is projective in [Fu2] is not necessary. it suffices to assume that $X$ is complete).
On the other hand, the statement was proved when $X$ is non-singular and $\Delta$ is a reduced divisor with only simple normal crossings in [Fu96] and when dim $X=2$ in [Fu3].
We note that, if $aL-(K_X+\Delta)$ is nef and big but not nef and log big on $(X,\Delta)$, there exists a counterexample due to Zariski (cf. [KMM, Remark 3-1-2]).

We shall prove the following result in this paper:

\proclaim{Main Theorem}
Assume that $(X,\Delta)$ is weakly Kawamata log terminal $($wklt$)$ and that every irreducible component of $\lfloor \Delta \rfloor$ is {\bf Q}-Cartier.
Let $L$ be a nef Cartier divisor such that $aL-(K_X+\Delta)$ is nef and log big on $(X,\Delta)$ for some $a \in$ {\bf N}.
Then Bs$\vert mL \vert = \emptyset$ for every $m \gg 0$.
\endproclaim

\flushpar
This implies a kind of "log abundance theorem":

\proclaim{Corollary}
If $(X,\Delta)$ is {\bf Q}-factorial weakly Kawamata log terminal and $K_X+\Delta$ is nef and log big on $(X,\Delta)$, then Bs$\vert mr(K_X+\Delta) \vert = \emptyset$ for every $m \gg 0$.
\endproclaim

\flushpar
We remark that, concerning the log abundance conjecture, the following facts are known:

\flushpar
(1) If $(X,\Delta)$ is klt and $K_X+\Delta$ is nef and big, then Bs$\vert mr(K_X+\Delta) \vert = \emptyset$ for every $m \gg 0$ (cf. [Ka84, Theorem 2.6]).

\flushpar
(2) If dim $X \leq 3$ and $K_X+\Delta$ is nef, then $K_X+\Delta$ is semi-ample ([Ka79], [Fuj], [Ka1], [Utah, 8.4], [KeMaMc]).

The theorem will be proved in Sect.2, by showing the following

\proclaim{Key Lemma}
If $(X,\Delta)$ is wklt and every irreducible component $S$ of $\lfloor \Delta \rfloor$ is {\bf Q}-Cartier, then $(S, ${\rm Diff}$(\Delta - S))$ is wklt and every irreducible component of $\lfloor $ {\rm Diff}$(\Delta - S) \rfloor$ is {\bf Q}-Cartier.
\endproclaim

\specialhead 1. Preliminaries
\endspecialhead

\flushpar
We collect some results that will be needed in the next section.

\proclaim{Proposition 1 (Shokurov's Connectedness Lemma)}
{\rm ([He, Lemma 2.2], \linebreak [Ka2, Theorem 1.4], cf. [Sh, 5.7], [Utah, 17.4])}
Let $W$ be a normal, complete algebraic variety over {\bf C}, $(W,\Gamma)$ a log variety which is log canonical and $g: V \to W$ a log resolution of $(W,\Gamma)$.
Then {\rm (}the support of the effective part of $\lfloor (g^*(K_W+\Gamma)-K_V) \rfloor${\rm )}$\cap g^{-1}(s)$ is connected for every $s \in W$.
\endproclaim

\proclaim{Proposition 2}
{\rm (cf. [Ka1, the proof of Lemma 3], [Ka84, Theorem 2.6])}
Assume that $(X,\Delta)$ is wklt.
Let $L$ be a nef Cartier divisor such that $aL-(K_X+\Delta)$ is nef and big for some $a \in$ {\bf N}.
If Bs$\vert mL \vert \cap \lfloor \Delta \rfloor = \emptyset$ for every $m \gg 0$, then Bs$\vert mL \vert = \emptyset$ for every $m \gg 0$.
\endproclaim

\proclaim{Proposition 3}
{\rm ([Sz])}
If $(X,\Delta)$ is divisorial log terminal, then $(X,\Delta)$ is wklt.
\endproclaim

\proclaim{Proposition 4 (Reid Type Vanishing)}
{\rm (cf. [Fu1], [Fu2, Proposition 1])}
\linebreak Assume that $(X,\Delta)$ is wklt.
Let $D$ be a {\bf Q}-Cartier integral Weil divisor.
If $D-(K_X+\Delta)$ is nef and log big on $(X,\Delta)$, then $H^i(X,{\Cal O}_X(D))=0$ for every 
$i>0$.
\endproclaim

\specialhead 2. Proof of Main Theorem
\endspecialhead

\flushpar
We simply sketch the proof of Main Theorem, since we proceed along the lines of that in [Fu2].

Let $S$ be an irreducible component of $\lfloor \Delta \rfloor$.
From [Utah, 17.5] (cf. [Sh, 3.8]), $S$ is normal.

Let $f:Y \to X$ be a log resolution of $(X,\Delta)$ such that the following conditions are satisfied:

(1) Exc($f$) consists of divisors,

(2) $K_Y+f_*^{-1}\Delta +F=f^*(K_X+\Delta)+E$,

(3) $E$ and $F$ are $f$-exceptional effective {\bf Q}-divisors such that Supp($E$) and Supp($F$) do not have common irreducible components,

(4) $\lfloor F \rfloor =0$.

\proclaim{Claim 1}
{\rm ([Sh, p.99])}
For any member $G \in$ {\bf Strata} $(f_*^{-1}\lfloor \Delta \rfloor)$, Exc$(f)$ does not include $G$.
\endproclaim

We put $S_0:=f_*^{-1}S$ and Diff$(\Delta - S):=(f \vert _{S_0})_*(f^*(K_X+\Delta)\vert _{S_0} -(K_Y+S_0)\vert _{S_0})$.

We note that $(K_X+\Delta)\vert _S = K_S +$ Diff$(\Delta - S)$.

\proclaim{Key Lemma}
$(S, ${\rm Diff}$(\Delta - S))$ is wklt and every irreducible component of $\lfloor $ {\rm Diff}$(\Delta - S) \rfloor$ is {\bf Q}-Cartier.
\endproclaim

\demo{Proof}
By the Subadjunction Lemma ([Sh, 3.2.2], cf. [KMM, Lemma 5-1-9]), {\rm Diff}$(\Delta - S) \ge 0$.
Here $\lfloor$ {\rm Diff} $(\Delta - S) \rfloor = (f \vert _{S_0})_*((f_*^{-1} \lfloor \Delta -S \rfloor )\vert _{S_0})$.
From [Sh, 3.2.3] and Proposition 3 or from [Ka1, Lemma 4], $(S, ${\rm Diff}$(\Delta - S))$ is wklt.

Let $D$ be an irreducible component of $\lfloor \Delta -S\rfloor$.
For $x \in D \cap S$, there exist $y_1 \in f_*^{-1}D$ and $y_2 \in S_0$ such that $f(y_1)=f(y_2)=x$.
Applying Proposition 1 to $(X,\{ \Delta \} +S+D)$ and $f$, we obtain $y_3 \in f_*^{-1}D \cap S_0$ such that $f(y_3)=x$.
Thus $(f \vert _{S_0})_*(f_*^{-1} D \vert _{S_0})=$ Supp $(D \vert _S)$ by Claim 1.

We put Diff$(\{ \Delta \} +D):=(f \vert _{S_0})_*(f^*(K_X+S+\{ \Delta \} +D)\vert _{S_0} -(K_Y+S_0)\vert _{S_0})$.

We note that $(K_X+S+\{ \Delta \} +D)\vert _S = K_S +$ Diff$(\{ \Delta \} +D)$.

Here Diff$(\{ \Delta \} +D) \quad \ge 0$ from [Sh, 3.2.2] (cf. [KMM, Lemma 5-1-9]) and \linebreak $\lfloor$ Diff$(\{ \Delta \} +D) \rfloor = \quad (f \vert _{S_0})_*(f_*^{-1} D \vert _{S_0})$.
Applying Proposition 1 to $(S,$ Diff$(\{ \Delta \} +D))$ and $f \vert_{S_0}$ , from the fact that every connected component of $f_*^{-1} D \vert _{S_0}$ is \linebreak irreducible, we deduce that every connected component of $(f \vert _{S_0})_*(f_*^{-1} D \vert _{S_0})$ is \linebreak irreducible.
As a result every irreducible component of $(f \vert _{S_0})_*(f_*^{-1} D \vert _{S_0})$ is {\bf Q}-Cartier, because $D \vert_S$ is {\bf Q}-Cartier.
\qed
\enddemo

\proclaim{Claim 2}
$aL \vert _S -(K_S+${\rm Diff}$(\Delta - S))$ is nef and log big on $(S, ${\rm Diff}$(\Delta - S))$.
\endproclaim

\demo{Proof}
The assertion follows from the fact that $\lfloor $ {\rm Diff}$(\Delta - S) \rfloor = (f \vert _{S_0})_*((f_*^{-1} \lfloor \Delta -S \rfloor )\vert _{S_0})$.
\qed
\enddemo

\proclaim{Claim 3}
$\left\vert mL \right\vert \big\vert_S = \big\vert mL \vert _S \big\vert$ for $m \ge a$.
\endproclaim

\demo{Proof}
We note that $mL-S$ is a {\bf Q}-Cartier integral divisor, $(X, \Delta -S)$ is wklt and $mL-S-(K_X+\Delta-S)$ is nef and log big on $(X, \Delta -S)$.
Thus $H^1(X,{\Cal O}_X(mL-S))=0$ from Proposition 4.
\qed
\enddemo

We complete the proof of the theorem by induction on dim$X$.

By Key Lemma and Claim 2 and by induction hypothesis, Bs$\big\vert mL \vert _S \big\vert = \emptyset$ for $m \gg 0$.
Thus, by Claim 3, Bs$\vert mL \vert \cap \lfloor \Delta \rfloor = \emptyset$ for every $m \gg 0$.
Consequently Proposition 2 implies the assertion.

\end
\Refs
\widestnumber\key{KeMaMc}

\ref\key Fuj
\by T. Fujita
\paper Fractionally logarithmic canonical rings of algebraic surfaces
\jour J. Fac. Sci. Univ. Tokyo Sect. IA Math. \vol 30 \pages 685--696
\yr 1984
\endref

\ref\key Fu96
\by S. Fukuda
\paper On base point free theorem
\jour Kodai Math. J. \vol 19 \pages 191--199
\yr 1996
\endref

\ref\key Fu1
\by S. Fukuda
\paper A generalization of the Kawamata-Viehweg vanishing theorem after Reid
\jour Commun. Algebra \vol 24 \pages 3265--3268
\yr 1996
\endref

\ref\key Fu2
\by S. Fukuda
\paper A base point free theorem of Reid type
\jour J. Math. Sci. Univ. Tokyo \vol 4 \pages 621--625
\yr 1997
\endref

\ref\key Fu3
\by S. Fukuda
\paper A base point free theorem for log canonical surfaces
\jour Osaka J. Math.
\toappear
\finalinfo (preprint alg-geom/9705005)
\endref

\ref\key He
\by S. Helmke
\paper On Fujita's conjecture
\jour Duke Math. J. \vol 88 \pages 201--216
\yr 1997
\endref

\ref\key Ka79
\by Y. Kawamata
\paper On the classification of non-complete algebraic surfaces
\jour Lect. Notes Math. \vol 732 \pages 215--232
\yr 1979
\endref

\ref\key Ka84
\by Y. Kawamata
\paper The cone of curves of algebraic varieties
\jour Ann. Math. \vol 119 \pages 603--633
\yr 1984
\endref

\ref\key Ka1
\by Y. Kawamata
\paper Log canonical models of algebraic 3-folds
\jour Internat. J. Math. \vol 3 \pages 351--357
\yr 1992
\endref

\ref\key Ka2
\by Y. Kawamata
\paper On Fujita's freeness conjecture for 3-folds and 4-folds
\jour Math. Ann. \vol 308 \pages 491--505
\yr 1997
\endref

\ref\key KMM
\by Y. Kawamata, K. Matsuda and K. Matsuki
\paper Introduction to the minimal model problem
\jour Adv. Stud. Pure Math. \vol 10 \pages 283--360
\yr 1987
\endref

\ref\key KeMaMc
\by S. Keel, K. Matsuki and J. McKernan
\paper Log abundance theorem for threefolds
\jour Duke Math. J. \vol 75 \pages 99--119
\yr 1994
\endref

\ref\key Utah
\by J. Koll\'ar et.al.
\paper Flips and abundance for algebraic threefolds
\jour Ast\'erisque \vol 211
\yr 1992
\endref

\ref\key Rd
\by M. Reid
\paper Commentary by M. Reid $($\S10 of Shokurov's paper "3-fold log-flips"$)$
\jour Russian Acad. Sci. Izv. Math. \vol 40 \pages 195--200
\yr 1993
\endref

\ref\key Sh
\by V. V. Shokurov
\paper 3-fold log-flips
\jour Russian Acad. Sci. Izv. Math. \vol 40 \pages 95--202
\yr 1993
\endref

\ref\key Sz
\by E. Szab\'o
\paper Divisorial log terminal singularities
\jour J. Math. Sci. Univ. Tokyo \vol 1 \pages 631--639
\yr 1994
\endref


\endRefs
